\documentclass[a4paper, 11pt,reqno]{amsart} 
\usepackage{amssymb}
\usepackage{amsmath}

\newtheorem{thm}{Theorem}[section]
\newtheorem{prop}{Proposition}[section]
\newtheorem{lem}{Lemma}[section]
\newtheorem{rem}{Remark}[section]
\newtheorem{cor}{Corollary}[section]

\makeatletter
   
          \@addtoreset{equation}{section}
\makeatother

\begin{document}

\baselineskip=1.3\baselineskip   

\title[Massless Dirac operators]{The zero modes and  zero resonances\\
of   massless Dirac operators }

\author{Yoshimi Sait\={o}
and
 Tomio Umeda\\}

\address{
Department of Mathematics \\
University of Alabama at Birmingham, 
Birmingham, AL 35294, USA}
\email{saito@math.uab.edu}

\address{Department of Mathematical Sciences,
University of Hyogo, Himeji 671-2201, Japan
}
\email{umeda@sci.u-hyogo.ac.jp}
\date{}

\maketitle

\vspace{20pt}

\textbf{Abstract.}  The zero modes  and  zero resonances of 
the Dirac operator $H=\alpha\cdot  D + Q(x)$ are discussed, where
 $\alpha= (\alpha_1, \, \alpha_2, \, \alpha_3)$ is
the triple of  $4 \times 4$ Dirac matrices,
$ D=\frac{1}{\, i \,} \nabla_x$, and 
$Q(x)=\big( q_{jk} (x) \big)$   is a $4\times 4$ Hermitian matrix-valued function
with
 $| q_{jk}(x) | \le C \langle x \rangle^{-\rho} $,  $\rho >1$.
 We shall show that every zero mode $f(x)$ is  continuous 
  on ${\mathbb R}^3$ and decays at infinity with the decay rate
    $|x|^{-2}$. 
   Also, we shall show that $H$
   has no zero resonance if $\rho > 3/2$.

 \vspace{15pt}


 \vspace{15pt}

\textbf{The 2000 Mathematical Subject Classification:} 35Q40, 35P99, 81Q10

\newpage

\section{Introduction}

This paper is concerned with the massless Dirac operator
\begin{equation} \label{eqn:1-1}
H= \alpha \cdot D + Q(x), \quad D=\frac{1}{\, i \,} \nabla_x,
\,\,\, x \in {\mathbb R}^3,
\end{equation}
where $\alpha= (\alpha_1, \, \alpha_2, \, \alpha_3)$ is
the triple of  $4 \times 4$ Dirac matrices
\begin{equation*}\label{eqn:1-2}
\alpha_j = 
\begin{pmatrix}
 \mathbf 0 &\sigma_j \\ \sigma_j &   \mathbf 0 
\end{pmatrix}  \qquad  (j = 1, \, 2, \, 3)
\end{equation*}
with  the $2\times 2$ zero matrix $\mathbf 0$ and
the triple of  $2 \times 2$ Pauli matrices
\begin{equation*}\label{eqn:1-3}
\sigma_1 =
\begin{pmatrix}
0&1 \\ 1& 0
\end{pmatrix}, \,\,\,
\sigma_2 =
\begin{pmatrix}
0& -i  \\ i&0
\end{pmatrix}, \,\,\,
\sigma_3 =
\begin{pmatrix}
1&0 \\ 0&-1
\end{pmatrix},
\end{equation*}
and $Q(x)$ is a $4\times 4$ Hermitian matrix-valued function
decaying at infinity.

We would like to emphasize that
one can regard the operator (\ref{eqn:1-1}) as 
a generalization of 
the operator
\begin{equation} \label{eqn:1-4}
\alpha\cdot \big(D - A(x) \big) + q(x) I_4 , 
\end{equation}
where $(q, A)$ is an electromagnetic potential and
$I_4$ is a $4\times 4$ identity matrix,
by  taking $Q(x)$ to be $- \alpha\cdot A(x) + q(x) I_4$.
In the case  where $q(x)\equiv 0$, the operator
(\ref{eqn:1-4}) becomes of the form
\begin{equation*}   \label{eqn:1-4-1}
\alpha\cdot \big(D - A(x) \big)
=
\begin{pmatrix}
 \mathbf 0 &\sigma \cdot (D -A(x)) \\ 
\sigma \cdot (D -A(x)) & \mathbf 0 
\end{pmatrix}.
\end{equation*}
The component $\sigma \cdot (D -A(x))$ is called 
the Weyl-Dirac operator. See Balinsky and Evans \cite{BalinEvan2}.

In the paper by Fr\"ohlich, Lieb 
and Loss \cite{FrohlichLiebLoss},
it was found that 
the existence of zero modes (i.e., eigenfunctions with
the zero eigenvalue) of the Weyl-Dirac operator plays a crucial role 
in the study of stability of Coulomb systems with magnetic
fields.
(For the precise definition of zero modes,
see Definition 1.1 in 
the latter part of this section.)
Loss and Yau \cite{LossYau} were 
the first to construct zero modes of the Weyl-Dirac 
operator and their results
were usefully applied in \cite{FrohlichLiebLoss}.
Since then,
the zero modes of
the Dirac operator
$\alpha
\cdot (D -A(x))$,
the Weyl-Dirac operator $\sigma \cdot (D -A(x))$
and the Pauli operator $\{ \sigma \cdot (D -A(x)) \}^2 + q(x) I_2$
have attracted 
a considerable attention.
It is now widely understood that
the zero modes have deep and fruitful
implications from the view point of
mathematics as well as physics. 
See Adam, Muratori and Nash \cite{AdamMuratoriNash1},
\cite{AdamMuratoriNash2}, \cite{AdamMuratoriNash3},
Balinsky and Evans \cite{BalinEvan1},
\cite{BalinEvan2},
\cite{BalinEvan3},
Elton \cite{Elton}
and,
Erd\"os and Solovej \cite{ErdosSolovej1}.
Also, see 
Bugliaro, Fefferman and Graf \cite{BugliaroFeffermanGraf}
and,
 Erd\"os and Solovej 
 \cite{ErdosSolovej2}, 
\cite{ErdosSolovej3},
where their main concern  
is Lieb-Thirring inequality
for the  Pauli operator with
a strong magnetic fields
and, as a by-product, 
an estimate of the density
of zero modes of the Weyl-Dirac operators
 was obtained.

As for the two-dimensional case, 
Aharonov and Casher \cite{AharonovCasher}
are believed to be the first to construct
examples of zero modes.
See
Erd\"os and Vougalter \cite{ErdosVougalter}, 
Rozenblum and Shirokov \cite{RozenblumShirokov}
and Persson \cite{Persson}
for recent works.

We should like to note that 
the operator (\ref{eqn:1-1}) also generalizes
the Dirac operator of the form
\begin{equation} \label{eqn:1-5}
\alpha\cdot D  + m(x) \beta  + q(x) I_4, 
\end{equation}
where $m(x)$ is considered to be a variable mass, and
$\beta$ is the $4\times 4$ matrix defined by
\begin{equation*} \label{eqn:1-6}
\beta =
\begin{pmatrix}
 I_2 &  \mathbf 0  \\ \mathbf 0 & - I_2 
\end{pmatrix}.
\end{equation*}
Spectral properties of the operator (\ref{eqn:1-5})
have been extensively studied in recent years.
See Kalf and Yamada \cite{KalfYamada},
Kalf, Okaji and Yamada \cite{KalfOkajiYamada},
Schmidt and Yamada \cite{SchmidtYamada},
Pladdy \cite{Pladdy} and
Yamada \cite{Yamada}.

Finally, we would like to emphasize 
the significant role of the zero modes 
and  zero resonances in the analysis of the asymptotic
behavior, around the origin of the
complex plane, of the  resolvent of the operator $H$ 
given by (\ref{eqn:1-1}).  
One can easily recognize the significance  as is
suggested by Jensen and Kato \cite{JensenKato}  on
the Schr\"odinger operator.

\vspace{10pt}

\noindent
\textbf{Notation.}

\noindent
The upper and lower half planes ${\mathbb C}_{\pm}$ are  
defined by
\begin{equation*}
{\mathbb C}_{+} 
  :=\{ \, z \in \mathbb C  \; | \;  \mbox{Im }z >0 \, \},
\qquad
{\mathbb C}_{-} 
  :=\{ \, z \in \mathbb C  \; | \;  \mbox{Im }z <0 \, \}
\end{equation*}
respectively.
By $S({\mathbb R}^3)$, we mean
 the Schwartz class of rapidly decreasing
functions on ${\mathbb R}^3$, and we set 
$\mathcal S = [S({\mathbb R}^3)]^4$.

By $L^2 =L^2({\mathbb R}^3)$, we mean the Hilbert space of
square-integrable functions on ${\mathbb R}^3$, and 
we introduce  a Hilbert space ${\mathcal L}^2$ by
     ${\mathcal L}^2 = [L^2({\mathbb R}^3)]^4$, where 
the inner product  
is given by
\begin{equation*}
    (f, g)_{{\mathcal L}^2}
 = \sum_{j=1}^4 (f_j, g_j)_{L^2}
\end{equation*}
for  $f = {}^t(f_1, f_2, f_3, f_4)$ 
and
$g = {}^t(g_1, g_2, g_3, g_4)$.
    
By $L^{2, s}({\mathbb R}^3)$, we mean the weighted $L^2$ space
defined by
\begin{equation*}
 L^{2, s}({\mathbb R}^3)
  :=\{ \, u \; | \; \langle x \rangle^s u 
 \in L^2({\mathbb R}^3) \, \}
\end{equation*}
with the inner product
\begin{equation*}
 (u, \, v)_{L^{2, s}}
:=
\int_{{\mathbb R}^3} 
   \langle x \rangle^{2s} u(x)  \, \overline{v(x)} \, dx, 
\end{equation*}
where
\begin{equation*}  \label{eqn:bracket-x}
\langle x \rangle = \sqrt{1 + |x|^2 \,}.
\end{equation*}
We introduce  the Hilbert space
     ${\mathcal L}^{2,s} = [L^{2,s}({\mathbb R}^3)]^4$ with
the inner product  
\begin{equation*}
    (f, g)_{{\mathcal L}^{2,s}}
 = \sum_{j=1}^4 (f_j, g_j)_{L^{2,s}}.
\end{equation*}
By $H^{\mu, s}({\mathbb R}^3)$, we mean the weighted 
Sobolev  space
defined by
\begin{equation*}
 H^{\mu, s}({\mathbb R}^3)
  :=\{ \, u \in {S^{\prime}({\mathbb R}^3)} \; | \; 
\langle x \rangle^s \langle D \rangle^{\mu}  \, u 
 \in L^2({\mathbb R}^3) \, \}
\end{equation*}
with the inner product
\begin{equation*}
 (u, \, v)_{H^{\mu, s}}
:=
\big(
   \langle x \rangle^{s} \langle D \rangle^{\mu} u ,
\langle x \rangle^{s}\langle D \rangle^{\mu} \, v
\big)_{L^2},
\end{equation*}
where 
\begin{equation}    \label{eqn:bracket-D}
\langle D \rangle = \sqrt{1 - \Delta \,}.
\end{equation}
In a similar fashion, 
we introduce  the Hilbert space
     ${\mathcal H}^{\mu,s} = [H^{\mu,s}({\mathbb R}^3)]^4$.
Note that $H^{\mu, 0}({\mathbb R}^3)$ coincides with
the Sobolev space of order $\mu\,$:$\,H^{\mu}({\mathbb R}^3)$,
and by  
${\mathcal H}^{\mu}$ we mean the
Hilbert space $[H^{\mu}({\mathbb R}^3)]^4$.
Also note
that
${\mathcal H}^{0,0} ={\mathcal L}^2$ and
${\mathcal H}^{0,s} ={\mathcal L}^{2,s}$.

By ${B}(\mu, s \, ; \,\nu, t )$, we mean the set of 
all bounded linear operators from
$H^{\mu, s}({\mathbb R}^3)$ into $H^{\nu, t}({\mathbb R}^3)$,
and by
$\mathcal{B}(\mu, s \, ; \,\nu, t)$,
the set of 
all bounded linear operators from
${\mathcal H}^{\mu, s}$ into 
${\mathcal H}^{\nu, t}$.
For an operator $W \in {B}(\mu, s \, ; \,\nu, t )$, 
we define a copy of $W \in \mathcal{B}(\mu, s \, ; \,\nu, t)$ by
\begin{gather}   \label{eqn:mapsto}
\begin{split}    
 {\mathcal H}^{\mu, s} \ni 
  f= \, &{}^t (f_1, \, f_2, \, f_3, \, f_4)    \\
{}&\mapsto
Wf ={}^t (Wf_1, \, Wf_2, \, Wf_3, \, Wf_4)
\in {\mathcal H}^{\nu, t}. 
\end{split}      
\end{gather}

\vspace{15pt}

\noindent
\textbf{Assumption (A).} 

\noindent
Each element $q_{jk}(x)$ 
($j, \, k =1, \, \cdots, \, 4$) of $Q(x)$ is 
a measurable function satisfying
\begin{equation} \label{eqn:2-1}
| q_{jk}(x) | \le C \langle x \rangle^{-\rho} 
\quad   ( \rho >1 ),
\end{equation}
where $C$ is a positive constant.
Moreover, $Q(x)$ is a Hermitian matrix for each 
$x \in {\mathbb R^3}$.

\vspace{15pt}

Note that, under Assumption (A), the Dirac operator (\ref{eqn:1-1})
is a self-adjoint operator in $\mathcal L^2$ with  
$\mbox{Dom}(H) = {\mathcal H}^1$.
The self-adjoint realization will be denoted 
by $H$ again. 
With an abuse of notation, we shall write $Hf$  
\textit{in the distributional sense} for
$f \in {\mathcal S}^{\prime}$ 
whenever it makes sense.

\vspace{15pt}

\noindent
D{\scriptsize EFINITION} 1.1. \ 
By a zero mode, we mean a function 
$f \in \mbox{Dom}(H)$ which satisfies
\begin{equation*}  \label{eqn:1-7}
Hf=0.
\end{equation*}
By a zero resonance, we mean a function 
$f \in {\mathcal L}^{2, -s}\setminus \mathcal L^2$, for
some $s > 0$,
 which satisfies
$\, Hf=0 \,$ in the distributional sense.

\vspace{15pt}

It is evident that a zero mode of $H$ is an eigenfunction of $H$
corresponding to the eigenvalue $0$, i.e.,
a zero mode is
 an element of $\mbox{Ker}(H)$,
 the kernel of the self-adjoint
operator $H$. 

It would seem that there is no decisive definition of
zero resonances.
A common understanding of zero resonances in the 
literature
is that a zero resonance is a non-$\mathcal L^2$
solution of $Hf =0$ in a space 
slightly larger than $\mathcal L^2$. 
(See, for example, Jensen and Kato \cite{JensenKato}.) 
In dealing with zero 
resonances in section 2 and later sections, 
we shall restrict ourselves to the case 
where $\rho >
3/2$ and 
$0 < s \le \min \{ 3/2, \, \rho- 1 \}$.

Balinsky and Evans \cite{BalinEvan2} is particularly interesting 
from our  view point in the sense that they dealt with 
the Weyl-Dirac  operator
$\sigma \cdot (D -A(x))$ and showed that the set of magnetic fields 
which give rise to zero modes
 is rather \lq\lq sparse.\rq\rq

In this paper, we investigate  the zero modes and
 zero resonances
of the operator $H$ in (\ref{eqn:1-1})
under Assumption(A).
Our goal  
is to establish a pointwise
estimate of the zero modes
as well as the continuity of the zero modes, and
also to show that the zero resonances do not
exist.

\vspace{15pt}

\section{Main results}

\begin{thm} \label{thm:th-2-2}
Suppose Assumption {\rm(A)} is satisfied. Let $f$ be  
a zero mode of the operator {\rm(\ref{eqn:1-1})}.
Then 

{\rm (i)} the  inequality 
\begin{equation}   \label{eqn:2-2}
|f(x)| \le C
\langle x \rangle^{-2}
\end{equation}
holds for all $x \in {\mathbb R}^3$, where the constant 
$C(=C_f)$ depends only on the zero mode $f$;

\vspace{4pt}

{\rm (ii)} the zero mode $f$ is a continuous function on ${\mathbb R}^3$.
\end{thm}
\vspace{15pt}

\begin{rem}
It is natural that zero modes
 exhibit only polynomial decays at infinity. 
In Loss and Yau \cite{LossYau},
they considered the Weyl-Dirac operator
$\sigma \cdot (D-A(x))$,
and  constructed two examples of  pairs of
 a vector potential $A$ and a zero mode $\psi$. 
One of their examples 
shows that $A(x) = O(|x|^{-2})$ 
and $\psi(x) = O(|x|^{-2})$ at infinity.
(Also, see examples in Adam, Muratori and
Nash \cite{AdamMuratoriNash1}.)
Thus,  it is true that the decay rate
 in Theorem \ref{thm:th-2-2}
is optimal at least for $\rho$ with  $1 < \rho \le 2$.
\end{rem}

\begin{rem}
In Bugliaro, Fefferman and Graf \cite{BugliaroFeffermanGraf} and,
 Erd\"os and Solovej \cite{ErdosSolovej2}, \cite{ErdosSolovej3},
 they 
established  estimates of the density of zero modes
of the Weyl-Dirac operator
$\sigma \cdot (D-A(x))$.
It is apparent that their estimates immediately imply
estimates of each zero mode.
These estimates of each zero
mode are, however, quite unclear in terms
of the decay rate at infinity 
because their estimates contain local 
lengthscales of the magnetic fields.
\end{rem}

\vspace{10pt}
Theorem \ref{thm:th-2-1} below means that zero
resonances do not exist
under the restriction on $s$ 
mentioned after Definition 1.1. 
Accordingly,
we need a larger $\rho$ than Theorem \ref{thm:th-2-2}.

\vspace{15pt}
\begin{thm} \label{thm:th-2-1}
Suppose Assumption {\rm(A)} is satisfied with
$\rho > 3/2$. If $f$ belongs to
${\mathcal L}^{2, \, -s}$ for some 
$s$ with $0 < s \le \min\{3/2, \, \rho - 1 \}$ and
satisfies $Hf=0$ in the 
distributional sense, then 
$f \in {\mathcal H}^1$.
\end{thm}

\vspace{15pt}

\section{A singular integral operator}  \label{sec:SIO}

One of the ingredients of the proofs of the main
theorems is a singular integral operator acting on 
four component vector functions.
The singular integral operator we deal with 
in this section is 
 defined by

\begin{equation} \label{eqn:sio-1}
A f(x) = 
\displaystyle{
\int_{{\mathbb R}^3}
  i \, \frac{\alpha \cdot (x - y)}{4\pi |x-y|^3} f(y) \, dy.
}
\end{equation}
for 
\begin{equation*}   \label{eqn:sio-2}
f = {}^t(f_1, \, f_2, \, f_3, \, f_4) 
\in {\mathcal L}^2,
\end{equation*}
where  $\alpha \cdot (x - y)$  means 
 the sum of the matrix operation $\alpha_j$
for the four-vector $(x_j - y_j) f$:
\begin{equation*}   \label{eqn:sio-2+}
\alpha \cdot (x - y) f  
 = \sum_{j=1}^3 \alpha_j (x_j - y_j) f. 
\end{equation*}

We shall need a few estimates of $A$ on
${\mathcal L}^2$ and on its subspaces.

\begin{lem}   \label{lem:siolem-1}
For each $f \in {\mathcal L}^2$,
 $Af(x)$ is defined
for a.e. $x \in {\mathbb R}^3$.
Moreover,  $A$ is a bounded operator from
${\mathcal L}^2$ to ${\mathcal L}^6$, i.e., there 
exists a constant  $C$ such that
\begin{equation*} \label{eqn:sio-3}
\Vert A f \Vert_{{\mathcal L}^6} 
\le
 C \Vert  f \Vert_{{\mathcal L}^2}
\end{equation*}
for all $f \in {\mathcal L}^2$.
\end{lem}

\begin{proof}
Since $\alpha_j$'s are 
  unitary matrices and satisfy the 
anti-commutation relation 
$\alpha_j \alpha_k  + \alpha_k \alpha_j = 2 \delta_{jk} I_4$,
 we have
\begin{equation*}    
| \alpha \cdot (x - y) \, f(y)|
= |x-y| \, |f(y)|.
\end{equation*}
(Note that $|x-y|$ and $|f(y)|$ 
are the Euclidean norms of ${\mathbb R}^3$ and  ${\mathbb R}^4$
respectively.)
Therefore we get
\begin{gather} 
\begin{split}     \label{eqn:sio-5}
|A f(x) | 
&\le 
\frac{1}{\,4\pi \,} 
\int_{{\mathbb R}^3}
\frac{1}{|x-y|^2}  |f(y)|\, dy  \\
{}
& = \frac{\,\pi \,}{2} I_1(\, |f| \, ) (x),
\end{split}
\end{gather}
where $I_1$ is the Riesz potential; see
Stein \cite[p.117]{Stein}. We shall appeal two well-known
facts
(Stein \cite[p.119]{Stein}) 
that  $I_1(\,  u \, ) (x)$
is finite for a.e. $x \in {\mathbb R}^3$ 
if $u \in L^2({\mathbb R}^3)$, 
and that
 $I_1$ is
a bounded operator from $L^2({\mathbb R}^3)$
to $L^6({\mathbb R}^3)$ 
(a special case of the Hardy-Littlewood-Sobolev inequality). 
These facts, together with (\ref{eqn:sio-5}),
yield the conclusions of the lemma.
\end{proof}

\vspace{10pt}

\begin{lem} \label{lem:siolem-2}
Let $s \ge 1$. Then

\begin{equation*} \label{eqn:yajima}
\Vert A f \Vert_{{\mathcal L}^2} 
\le
 C \Vert  f \Vert_{{\mathcal L}^{2, s}}
\end{equation*}
for all $f \in {\mathcal L}^{2, s}$.
\end{lem}

\begin{proof}
In view of (\ref{eqn:sio-5}), 
it is sufficient to show that
the
Riesz potential
$I_1$ is a bounded operator from 
$L^{2, s}({\mathbb R}^3)$ to $L^{2}({\mathbb R}^3)$.

Let $u \in S({\mathbb R}^3)$. Then we have 
\begin{equation}  \label{eqn:yajima-1}
I_1(u) =\overline{\mathcal F}
\, [ \frac{1}{\, 2\pi |\xi| \,} ]  
{\mathcal F} u,
\end{equation}
where ${\mathcal F}$ and $\overline{\mathcal F}$ denote
the Fourier transform and the inverse Fourier 
transform respectively (see \cite[p.117]{Stein}), and 
\begin{equation*}   \label{eqn:yajima-2}
{\mathcal F} u (\xi) = \hat u (\xi)
= \frac{1}{(2\pi)^{3/2}}
\int_{{\mathbb R}^3} e^{-ix\cdot \xi} \, u(x) \, dx.
\end{equation*}
It follows from (\ref{eqn:yajima-1}) and the
Plancherel theorem that
\begin{equation}   \label{eqn:yajima-3}
\Vert I_1(u) \Vert_{L^2}^2 
=
 \frac{1}{4\pi^2}
\int_{{\mathbb R}^3}   \frac{1}{\,|\xi|^2 \,}    \,
   |\hat u (\xi)|^2
 \, d\xi.
\end{equation}
If we apply the Hardy inequality (which is referred to as 
the uncertainty
principle lemma  in \cite[p.169]{ReedSimon};
also see \cite[p.4.50]{Kuroda-A}) to the
right hand side of  (\ref{eqn:yajima-3}), we get
\begin{gather} 
\begin{split}     \label{eqn:yajima-4}
\frac{1}{\, 4\pi^2 \,}
\int_{{\mathbb R}^3}   \frac{1}{\,|\xi|^2 \,}    \,
   |\hat u (\xi)|^2
 \, d\xi
& \le
\frac{1}{\, \pi^2 \,}
\int_{{\mathbb R}^3}   
   | \nabla_{\!\!\xi} \, \hat u (\xi)|^2   \, d\xi    \\
\noalign{\vskip 4pt}
&=
\frac{1}{\, \pi^2 \,}
\int_{{\mathbb R}^3}   
   | x  \, u (x)|^2   \, dx. 
\end{split}
\end{gather}
Combining (\ref{eqn:yajima-3}) and (\ref{eqn:yajima-4}),
we obtain
\begin{equation}   \label{eqn:yajima-5}
\Vert I_1(u) \Vert_{L^2} 
\le 
 \frac{1}{\, \pi \, } \Vert x  u \Vert_{L^2} 
\le 
 \frac{1}{\, \pi \, } 
\Vert   u    \Vert_{L^{2, s}} 
\end{equation}
for $u \in S({\mathbb R}^3)$,
where we have used the hypothesis $s \ge 1$.
Since $S({\mathbb R}^3)$ is dense in 
$L^{2, s}({\mathbb R}^3)$, it follows from
(\ref{eqn:yajima-5}) that
$I_1$ is bounded from 
$L^{2, s}({\mathbb R}^3)$ to $L^{2}({\mathbb R}^3)$.
\end{proof}

\vspace{10pt}

We introduce a class of functions which
is necessary to establish an 
${\mathcal L}^{\infty}$ estimate 
of the operator $A$. For $q\ge 1$, 
we define
\begin{equation*}  \label{eqn:sio-6}
L^q_{ul}(\mathbb R^3)
=
\big\{\,
u \in L^q_{loc}(\mathbb R^3) \; \big| \;
 \Vert u \Vert_{L^q_{ul}} := 
\sup_{x \in {\mathbb R}^3}  
\Vert u \Vert_{L^q{( B(x; \, 1) )}} < \infty
\, \big\},
\end{equation*}
where $B(x;  1) 
= \{ \, y \in \mathbb R^3 \; | \;\, |x-y| \le 1 \, \}$, and
define
\begin{equation*}  \label{eqn:sio-7}
{\mathcal L}^q_{ul}
=
 \big[ L^q_{ul}(\mathbb R^3) \big]^4,
\quad 
\Vert  f \Vert_{{\mathcal L}^q_{ul}}
= \sum_{k=1}^4 
\Vert f_k \Vert_{L^q_{ul}}.
\end{equation*}

\vspace{10pt}

\begin{lem} \label{lem:siolem-3}
Let $1 <p < 3 < q < +\infty$. Then there exists a constant 
$C_{pq}$ such that
\begin{equation*}  \label{eqn:sio-8}
\Vert A f \Vert_{{\mathcal L}^{\infty}}
\le C_{pq}
 \big( 
\Vert  f \Vert_{{\mathcal L}^p}
+\Vert  f \Vert_{{\mathcal L}^q_{ul}}
\big).
\end{equation*}
for all  $f \in {\mathcal L}^p \cap {\mathcal L}^q_{ul}$.
In particular,
\begin{equation*}  \label{eqn:sio-8}
\Vert A f \Vert_{{\mathcal L}^{\infty}}
\le C_{pq}
 \big( 
\Vert  f \Vert_{{\mathcal L}^p}
+\Vert  f \Vert_{{\mathcal L}^q}
\big).
\end{equation*}
for all  $f \in {\mathcal L}^p \cap {\mathcal L}^q$.
\end{lem}
\begin{proof}
By virtue of (\ref{eqn:sio-5}), we only have to prove 
that there exists a constant $C_{pq}^{\prime}$ such that
\begin{equation}   \label{eqn:sio-9}
\Vert I_1(u) \Vert_{{L}^{\infty}}
\le C_{pq}^{\prime}
 \big( 
\Vert  u \Vert_{{L}^p}
+\Vert u \Vert_{{L}^q_{ul}}
\big)  
\;\;\mbox{ for }  u \in 
L^p(\mathbb R^3) \cap L^q_{ul}(\mathbb R^3).
\end{equation}
Since each $u \in 
L^p(\mathbb R^3) \cap L^q_{ul}(\mathbb R^3)$ can be
decomposed as
\begin{align*}  
{}&u =v_+ -v_- + i (w_+ -w_-),  \\
{}&v_{\pm} \ge 0, \;\; w_{\pm}\ge 0, 
\;\;
 v_{\pm},\,  w_{\pm} 
\in  L^p(\mathbb R^3) \cap L^q_{ul}(\mathbb R^3).
\end{align*}
we shall prove (\ref{eqn:sio-9}) for $u\ge0$.

Let $u \in 
L^p(\mathbb R^3) \cap L^q_{ul}(\mathbb R^3)$ be given,
and let satisfy $u\ge0$. 
Then one can find a sequence $\{ \, \varphi_n \, \} \subset
C_0^{\infty}(\mathbb R^3)$ such that
\begin{equation}   \label{eqn:sio-13}
0 \le \varphi_n \le u, \quad  \varphi_n  \to u  \;\; 
\mbox{ in } L^p(\mathbb R^3).
\end{equation}
(First cut $u$ as ${\chi}_{B(0;\, n)}(x) u$ by multiplying
a characteristic function ${\chi}_{B(0;\, n)}$ 
of the ball $B(0;\, n)$ with center at 
the origin and radius $n$, then use
the mollifier.) 
For each $n$, 
we decompose as
\begin{gather} 
\begin{split}    \label{eqn:sio-14}
I_1(\varphi_n)(x)
&= 
\int_{\mathbb R^3} \frac{1}{\, 2\pi^2 |x-y|^2 \,} \,
  \varphi_n (y) \, dy        \\
\noalign{\vskip 4pt}
&= h_0 * \varphi_n (x)  
  +  h_1 * \varphi_n (x),     
\end{split}
\end{gather}
where
\begin{gather*}
\begin{split}   
h_0(x) &= \chi_{B(0;\, 1)}(x)
 \frac{1}{\,  2\pi^2 |x|^2 \,} ,      \\
h_1(x) &= (1 -\chi_{B(0;\, 1)}(x) )
 \frac{1}{\,  2\pi^2 |x|^2 \,}.
\end{split}
\end{gather*}
(One should note that the integral on the right hand side of
(\ref{eqn:sio-14})  converges because of the fact that
$\varphi_n \in C_0^{\infty}(\mathbb R^3)$.)
If we apply the H\"older inequality to $h_0 * \varphi_n$,
then we get
\begin{gather}   
\begin{split}\label{eqn:sio-17}
|h_0 * \varphi_n(x) | 
&\le
 \frac{1}{\,  2\pi^2 \,}  
\Big\{
 \int_{|x-y|\le 1} \frac{1}{\,  |x-y|^{2  q^{\prime}} \,}  \, dy
\Big\}^{1/  q^{\prime}} \,
\Vert \varphi_n \Vert_{L^q(B(x;1))}   \\
\noalign{\vskip 4pt}
{}&
\quad \le  
C_{q}^{\prime}  \Vert u \Vert_{L^q(B(x;1))}    \;\;\quad
\big( 
   \, \frac{1}{ q^{\prime}}= 1 - \frac{1}{q}   \,
     \big),
\end{split}
\end{gather}
where we have used the fact that
 $2  q^{\prime} < 3 \;\; (\because q >3$ by assumption) 
and  (\ref{eqn:sio-13}).
Similarly, if we apply the H\"older inequality
to $h_1 * \varphi_n$, we obtain
\begin{gather}   
\begin{split}\label{eqn:sio-18}
|h_1 * \varphi_n(x) | 
&\le
 \frac{1}{\,  2\pi^2 \,}  
\Big\{
 \int_{|x-y|\ge 1} \frac{1}{\,  |x-y|^{2  p^{\prime}} \,}  \, dy
\Big\}^{1/ p^{\prime}} \,
\Vert \varphi_n \Vert_{L^p({\mathbb R})}   \\
\noalign{\vskip 4pt}
{}&
\quad \le  
C_p^{\prime}  \Vert u \Vert_{L^p({\mathbb R}^3)}    \;\;\quad
\big( 
   \, \frac{1}{ p^{\prime}}= 1 - \frac{1}{p}   \,
     \big),
\end{split}
\end{gather}
where we have used the fact that
 $2  p^{\prime} > 3 \;\; (\because 1< p < 3$ by assumption) 
and  (\ref{eqn:sio-13}).
If follows from (\ref{eqn:sio-14}), 
(\ref{eqn:sio-17})  and (\ref{eqn:sio-18})  that
\begin{gather}   
\begin{split}\label{eqn:sio-19}
|I_1(\varphi_n)(x)|
& \le  
C_{pq}^{\prime} 
\big(  \Vert u \Vert_{L^q(B(x;1))}  
 + \Vert u \Vert_{L^p({\mathbb R}^3)}  \, \big) \\
& \le
C_{pq}^{\prime}
\big(  \Vert u \Vert_{L^q_{ul}}  
 + \Vert u \Vert_{L^p}   \big)
\end{split}
\end{gather}
Recall that  the Riesz potential
$I_1$ is a bounded operator 
from $L^p({\mathbb R}^3)$ to $L^r({\mathbb R}^3)$,
  ($r^{-1} = p^{-1} - 3^{-1}$), 
because of the Hardy-Littlewood-Sobolev inequality.
This fact, together with (\ref{eqn:sio-13}), implies that
there exists a subsequence $\{ \varphi_{n^{\prime}} \}$
such that
$I_1(\varphi_{n^{\prime}}) (x) \to I_1(u)(x)$ for a.e.  
$x \in {\mathbb R}^3$. 
Thus taking the limit of (\ref{eqn:sio-19}),
along with the subsequence,
gives (\ref{eqn:sio-9}).
\end{proof}

\vspace{15pt}

\section{Estimates of the resolvents}  \label{sec:resolvent}

Another  ingredient of the proofs of the main theorems
is the limiting absorption principle (LAP) for
the free Dirac operator
\begin{equation} \label{eqn:reslvnDrc-00}
H_0= \alpha \cdot D.
\end{equation}
We note that $H_0$ with $\mbox{Dom}(H_0) ={\mathcal H}^1$ is 
a self-adjoint operator in ${\mathcal L}^2$. 
The self-adjoint realization will
be denoted by $H_0$ again. It is well-known
that the spectrum $\sigma(H_0)$ coincides with
 the whole 
real line $\mathbb R$.
With an abuse of notation again, we shall write
$H_0f$ for $f \in {\mathcal S}^{\prime}$.

We first prepare a lemma, which will be needed in 
the proof of Theorem \ref{thm:th-2-1} in section \ref{sec:pfthm}.

\begin{lem} \label{lem:sobolev-1}
If $f \in {\mathcal L}^{2}$ and  
 $(\alpha \cdot D) f \in {\mathcal L}^{2}$,
 then $f \in{\mathcal H}^{1}$.
\end{lem}
\begin{proof}
We take the Fourier transform of $(\alpha \cdot D)f$,
and we have
\begin{equation}    \label{eqn:proof-5}
{\mathcal F}[ (\alpha \cdot D)f]
= (\alpha \cdot \xi) \hat f.
\end{equation}
Then by using assumption of the lemma 
and (\ref{eqn:proof-5}), we see
that
\begin{gather}   \label{eqn:proof-6}
\begin{split}
+ \infty > \Vert (\alpha \cdot D)f \Vert_{{\mathcal L}^2}^2
&=
\int_{{\mathbb R}^3}  
| (\alpha \cdot \xi) \hat f  (\xi) | ^2 \, d\xi   \\
&=
\int_{{\mathbb R}^3}
 < \!\! (\alpha \cdot \xi) \hat f  (\xi) ,
 (\alpha \cdot \xi) \hat f  (\xi) \!\! >_{\!{}_{\mathbb C}}  \, d\xi \\
&=
\int_{{\mathbb R}^3}
 |\xi|^2 |\hat f  (\xi)|^2  \, d\xi,
\end{split}
\end{gather}
where $< \cdot\, ,  \cdot>_{\!{}_{\mathbb C}}$ 
denotes the inner product of ${\mathbb C}^4$.
In the third equality of (\ref{eqn:proof-6}), we have used the
anti-commutation  relation.
Since, by assumption of 
the lemma, $f \in {\mathcal L}^{2}$,
 the conclusion of the lemma
follows  from (\ref{eqn:proof-6}).
\end{proof}

The task in the rest of this section is to
prove the following theorem,
which is essential in the proofs
of the main theorems in section \ref{sec:pfthm}.

\begin{thm} \label{thm:riesz3-5}
If $f \in {\mathcal L}^{2, -3/2}$ 
  and  $H_0 f \in {\mathcal L}^{2, s}$ 
for some $\, s > 1/2$,
 then $AH_0 f = f$.
\end{thm}

As was indicated at the beginning of this section,
the ingredient of the proof of
Theorem \ref{thm:riesz3-5} is
the
LAP for the free Dirac operator $H_0$.
 Our idea of proving it 
is based on a decomposition  of the
resolvent 
\begin{equation}   \label{eqn:reslvnDrc-0}
R_0(z) = (H_0 +z ) \varGamma_0(z^2)
  \;\; \mbox{ on } \; {\mathcal L}^2, 
\quad \mbox{Im}\, z \not=0,
\end{equation}
where
\begin{equation} \label{eqn:reslvnDrc-1}
R_0(z) = (H_0 - z )^{-1},
\end{equation}
and  $\varGamma_0(z)$ in (\ref{eqn:reslvnDrc-0})
denotes the copy of 
the resolvent $\varGamma_0(z) =(-\Delta -z)^{-1}$ of
the negative Laplacian
\begin{equation*}   \label{eqn:laplace-0}
-\Delta = -\left(
\frac{{\partial}^2}{\partial x_1^2} 
+
\frac{{\partial}^2}{\partial x_2^2} 
+
\frac{{\partial}^2}{\partial x_3^2} 
\right).
\end{equation*}
See  (\ref{eqn:mapsto}) for the definition
of the copy of an operator.
In other words,
we shall not distinguish between
$\varGamma_0 (z)$ in 
$L^2({\mathbb R}^3)$ 
and 
$\varGamma_0 (z)$ in ${\mathcal L}^2$.
We believe this will not cause any confusion.

A formal computation, based on
the anti-commutation relation, shows that
\begin{equation}   \label{eqn:diracsquare}
H_0^2 = -\Delta I_4,
\end{equation}
from which one can deduce 
(\ref{eqn:reslvnDrc-0}).
The decomposition (\ref{eqn:reslvnDrc-0}) was
first exploited in 
Balslev and Hellfer \cite{BalslevHellfer}.
Similar decomposition was also adopted in Pladdy, Sait\={o} and
Umeda \cite{PSU1}, \cite{PSU2}.

We shall divide the rest of this section
into two subsections, 
because the proof of Theorem \ref{thm:riesz3-5}
is lengthy.

\vspace{10pt}
\subsection{The resolvent of the negative Laplacian}

In this subsection, we shall state several lemmas,
which are actually well-known and 
reproductions  of results in
Jensen and Kato \cite{JensenKato} and 
Kuroda \cite{Kuroda-A}, \cite{Kuroda-I}.
We shall do this for our later purpose
as well as for the reader's convenience.

We first recall that the resolvent of $-\Delta$ can be 
represented as an integral operator:
\begin{equation}  \label{eqn:resolventLapc}
\varGamma_0 (z) u(x) =
\int_{{\mathbb R}^3} 
\frac{e^{i \sqrt{z} |x-y|} }{\, 4\pi |x-y| \,} u(y) \, dy,
\qquad
u \in L^2({\mathbb R}^3 )
\end{equation}
for $z \in \mathbb C \setminus [0, \,+\infty)$, 
where
$\mbox{Im} \sqrt{z} >0$.

We next recall  well-known inequalities 
(e.g., \cite[Appendix A]{Eckardt}, 
\cite[p.162]{Kuroda-I},
\cite[Lemma 11.1]{Umeda}),
 which will be repeatedly used in the present
paper: 
\begin{equation}   \label{eqn:ekku}
\int_{{\mathbb R}^3}
\frac{1}{ \, |x - y|^2  \, 
 \langle y \rangle^{\gamma}\,}  \, dy 
\le
C_{\gamma} 
\begin{cases}  
\langle x \rangle^{-\gamma + 1} 
& \mbox{if } \; 1  < \gamma < 3, 
\\ {}    & \\
\langle x \rangle^{-2} \log ( 1 + \langle x \rangle )
  & \mbox{if } \; \gamma = 3,
\\ {}    & \\
\langle x \rangle^{-2}
  & \mbox{if } \;  \gamma > 3.
\end{cases}
\end{equation}

\vspace{10pt}

\begin{lem}   \label{lem:hs}
Let  $s$, $s^{\prime} > 1/2$ and $s+ s^{\prime} >2$.
Then 
\begin{equation}  \label{eqn:hs-2}
\iint_{  {\mathbb R}^3 \times {\mathbb R}^3 } 
\,  \langle x \rangle^{-2s^{\prime}} 
     \frac{1}{ \, |x - y|^2 \,}
  \langle y \rangle^{-2s}  \, dxdy < +\infty.
\end{equation}
\end{lem}

\begin{proof}
We may assume, with no loss of generality, that $s < 3$. 
Then, by an inequality in (\ref{eqn:ekku}), we have
\begin{equation}  \label{eqn:hs-3}
\int_{ {\mathbb R}^3} 
     \frac{1}{ \, |x - y|^2 \,}
  \langle y \rangle^{-2s}  \, dy 
\le 
C_s \langle x \rangle^{-2s +1}.
\end{equation}
Since $-2s^{\prime} -2s +1 < -3$ by assumption of 
the lemma, we see that (\ref{eqn:hs-3})
implies   (\ref{eqn:hs-2}).
\end{proof}

\vspace{10pt}

It follows from  (\ref{eqn:resolventLapc}) and
 Lemma \ref{lem:hs} that 
the operator
\begin{equation}   \label{eqn:khs-0}
K(z) :=\langle x \rangle^{-s^{\prime}} \!
 \varGamma_0(z)  
\langle x \rangle^{-s},
\end{equation}
which is represented as
\begin{equation*}   \label{eqn:khs-1}
K(z) u (x) =
\int_{{\mathbb R}^3} 
\langle x \rangle^{-s^{\prime}}
\frac{e^{i \sqrt{z} |x-y|} }{\, 4\pi |x-y| \,}
\langle y \rangle^{-s} u(y) \, dy,
\end{equation*}
belongs to the Hilbert-Schmidt class  on
${L}^2({\mathbb R}^3)$ for 
$z \in \mathbb C \setminus [0, \,+\infty)$:
\begin{equation*}   \label{eqn:khs-1+}
\Vert 
K(z)        
  \Vert^2_{\mbox{\tiny{HS}}} 
=
\iint_{ {\mathbb R}^3 \times {\mathbb R}^3 } 
\,  \langle x \rangle^{-2s^{\prime}} 
\Big| \,
\frac{e^{i \sqrt{z} |x-y|} }{\, 4\pi |x-y| \,} 
 \, \Big|^2
  \langle y \rangle^{-2s}  \, dxdy < +\infty,
\end{equation*}
where $\Vert \cdot \Vert_{\mbox{\tiny{HS}}}$ denotes
the Hilbert-Schmidt norm.
Note that
\begin{gather}  \label{eqn:hs-4}
\begin{split}
\Vert 
K(z_1)  - K(z_2)        
  \Vert^2_{\mbox{\tiny{HS}}}   
       \qquad\qquad\qquad\quad \qquad\qquad
         \qquad\qquad\qquad\quad  \\
\noalign{\vskip 4pt}
=
\iint_{  {\mathbb R}^3 \times {\mathbb R}^3 } 
\,  \langle x \rangle^{-2s^{\prime}} 
\Big| \,
\frac{e^{i \sqrt{z_1} |x-y|} }{\, 4\pi |x-y| \,} 
-     
\frac{e^{i \sqrt{z_2} |x-y|} }{\, 4\pi |x-y| \,}
 \, \Big|^2
  \langle y \rangle^{-2s}  \, dxdy
\end{split}
\end{gather}
for all 
$z_1$, $z_2 \in \mathbb C \setminus [0, \,+\infty)$.
It follows from (\ref{eqn:hs-4}) that
$K(z)$
is continuous,
 with respect to the Hilbert-Schmidt
norm topology, on 
$\mathbb C \setminus [0, \,+\infty)$.
Furthermore, we can deduce 
from  (\ref{eqn:hs-2}), 
(\ref{eqn:hs-4})  and 
Lebesgue's  convergence theorem 
that
$K(z)$ can be continuously extended, 
with respect to the Hilbert-Schmidt
norm topology, as follows:
\begin{equation}   \label{eqn:khs-2}
\widetilde K(z) 
=
\begin{cases}  
K(z) 
& \mbox{if } \; z \in \mathbb C \setminus [0, \, +\infty), 
\\ {}    & \\
K^+(\lambda) 
  & \mbox{if } \;  z= \lambda + i0, \; \lambda \ge 0,
\\ {}    & \\
K^-(\lambda) 
  & \mbox{if } \; z= \lambda - i0, \; \lambda \ge 0,  
\end{cases}
\end{equation}
where $K^+(\lambda)$ and $K^-(\lambda)$ for 
$\lambda >0$
are defined by
\begin{equation}  \label{eqn:khs-3}
K^{\pm} (\lambda) u(x) =
\int_{{\mathbb R}^3} 
\langle x \rangle^{-s^{\prime}} \,
\frac{\, e^{\pm i \sqrt{\lambda} |x-y|} \, }{4\pi |x-y|} 
 \, \langle y \rangle^{-s} 
u(y) \, dy,
\end{equation}
and 
\begin{equation}  \label{eqn:khs-4}
K^+(0)u(x)= K^-(0)u(x)
=\int_{{\mathbb R}^3} 
 \langle x \rangle^{-s^{\prime}} 
\frac{1}{\, 4\pi |x-y| \,}  
 \langle y \rangle^{-s}  u(y) \, dy.
\end{equation}

For a later purpose, it is convenient to introduce 
a subset of the Riemann surface of $\sqrt{z \,}$ 
as follows:
\begin{gather}   \label{eqn:riemann-1}
\begin{split}
{}&\;\; \Pi_{(0, \, +\infty)}     \\
&:=
\big( \mathbb C \setminus (0, \, +\infty) \big)
\cup
\{ \,  z = \lambda + i0 \; | \; \lambda >  0  \, \}
\cup
\{ \, z = \lambda - i0 \; | \; \lambda >  0  \, \}.
\end{split}
\end{gather}
Thus, we can say that $\widetilde K(z)$ defined 
by (\ref{eqn:khs-2}) -- (\ref{eqn:khs-4}) 
is continuous on $\Pi_{(0, \, +\infty)}$
with respect to
the Hilbert-Schmidt norm topology.

In view of (\ref{eqn:khs-0}), 
we see that $\varGamma_0(z)$, 
$z \in \mathbb C \setminus [0, \,+\infty)$, 
is a Hilbert-Schmidt operator from
 $L^{2, s}({\mathbb R}^3)$ to 
$L^{2, -s^{\prime}}({\mathbb R}^3)$.
Hence, in particular,
$\varGamma_0(z) \in {B}(0, s \, ; \,0 , -s^{\prime})$, 
and
\begin{equation*}   \label{eqn:hs-5-0}
\Vert
\varGamma_0(z)\Vert_{{B}(0, s \, ; \,0 , -s^{\prime})}   
\le
\Vert
K(z) 
\Vert_{\mbox{\tiny{HS}}}  .
\end{equation*}
Since we have the inequality
\begin{gather}
\begin{split}   \label{eqn:hs-5}
\Vert
\varGamma_0(z_1)  - \varGamma_0(z_2)
\Vert_{{B}(0, s \, ; \,0 , -s^{\prime})}   
\le
\Vert
K(z_1)  - K(z_2) 
\Vert_{\mbox{\tiny{HS} } }  \\ 
\noalign{\vskip 4pt}
(z_1, \; z_2 \in   \mathbb C \setminus [0, \, &+\infty)), 
\end{split}
\end{gather}
we conclude from  (\ref{eqn:khs-2}) 
  and (\ref{eqn:hs-5}) that 
$\varGamma_0(z)\in 
{B}(0, s \, ; \,0 , - s^{\prime})$ 
can be continuously extended
 as
follows:
\begin{equation}   \label{eqn:laplace-1}
\widetilde\varGamma_0(z) 
=
\begin{cases}  
\varGamma_0(z) 
& \mbox{if } \; z \in \mathbb C \setminus [0, \, +\infty), 
\\ {}    & \\
\varGamma_0^+(\lambda) 
  & \mbox{if } \;  z= \lambda + i0, \; \lambda \ge 0,
\\ {}    & \\
\varGamma_0^-(\lambda) 
  & \mbox{if } \; z= \lambda - i0, \; \lambda \ge 0,  
\end{cases}
\end{equation}
where
\begin{equation*}   \label{eqn:laplace-2}
\varGamma_0^{\pm}(\lambda) :=
\langle x \rangle^{s^{\prime}} 
 K^{\pm}(\lambda)
\langle x \rangle^{s}
=
\lim_{\epsilon \downarrow 0}   
\varGamma_0(\lambda \pm i \epsilon)  
 \; \mbox{ in } \;   
{B}(0, s \, ; \,0 , - s^{\prime}).
 \end{equation*}
We must remark that 
\begin{equation}  \label{eqn:resolventLapc-1}
\varGamma_0^+(0)u(x)= \varGamma_0^-(0)u(x)
=\int_{{\mathbb R}^3} 
\frac{1}{\, 4\pi |x-y| \,} u(y) \, dy,
\end{equation}
and that
\begin{equation*}  \label{eqn:resolventLapc-2}
\varGamma_0^{\pm} (\lambda) u(x) =
\int_{{\mathbb R}^3} 
\frac{e^{\pm i \sqrt{\lambda} |x-y|} }{\, 4\pi |x-y| \,} u(y) \, dy.
\end{equation*}
Thus
\begin{equation*}  \label{eqn:resolventLapc-3}
\widetilde\varGamma_0(\lambda + i0)
 \not= \widetilde\varGamma_0(\lambda -i0), \quad \lambda >0.
\end{equation*}
Note that the equality (\ref{eqn:resolventLapc-1}) allows us 
to use the notation
\begin{equation}  \label{eqn:resolventLapc-3-1}
\widetilde\varGamma_0(0)
\big( =\varGamma_0^+(0)= \varGamma_0^-(0) \big).
\end{equation}
With the notation introduced in (\ref{eqn:riemann-1}),
we can  say that
$\widetilde\varGamma_0(z)$ is a\linebreak  
${B}(0, s \, ; \, 0, - s^{\prime})$-valued
continuous function on $\Pi_{(0, \, +\infty)}$.

The following lemmas \ref{lem:lap-1} and 
\ref{lem:lap-1-1}  
are variants of Lemma 2.1 of  
Jensen and Kato \cite{JensenKato},
although we shall give their proofs. 

\begin{lem} \label{lem:lap-1}
Let $s$, $s^{\prime}$ satisfy
the same assumption as in Lemma \ref{lem:hs},
 and let $\mu \in \mathbb R$.
 Then $\widetilde\varGamma_0(z)$ is a  
${B}(\mu, s \, ; \,\mu, - s^{\prime})$-valued
continuous function
on $\Pi_{(0, \, +\infty)}$.
\end{lem}

\begin{proof} 
As was mentioned before,
$\widetilde\varGamma_0(z)$ defined by (\ref{eqn:laplace-1})
is 
a ${B}(0, s \, ; \,0 ,-s^{\prime})$-valued
continuous function  on $\Pi_{(0, \, +\infty)}$.
 Then the lemma directly
follows from the inequalities
\begin{gather}   \label{eqn:laplace-3}
\begin{split}
\Vert
\varGamma_0(z)
\Vert_{{B}(\mu, s \, ; \,\mu , -s^{\prime})}
\le
\Vert
\varGamma_0(z)
\Vert_{{B}(0, s \, ; \,0 , -s^{\prime})},
\quad   
z \in   \mathbb C
 \setminus [0, \,+\infty)   
\end{split}
\end{gather}
and
\begin{gather}    \label{eqn:laplace-3+}
\begin{split}
\Vert
\varGamma_0(z_1) - \varGamma_0(z_2)
\Vert_{{B}(\mu, s \, ; \,\mu , -s^{\prime})}
\le
\Vert
\varGamma_0(z_1) - \varGamma_0(z_2)
\Vert_{{B}(0, s \, ; \,0 , -s^{\prime})}     \\
\noalign{\vskip 4pt}   
( z_1, \, z_2 \in 
\mathbb C \setminus [0, \,+\infty)).    \qquad\qquad\qquad
\end{split}
\end{gather}
In order to show  (\ref{eqn:laplace-3}),
we shall use the fact
that
\begin{equation}  \label{eqn:laplace-4}
\langle D \rangle^{\mu} \varGamma_0(z) u
=  \varGamma_0(z) \langle D \rangle^{\mu} u
\end{equation}
for $u\in  S({\mathbb R}^3)$ and 
$z\in \mathbb C \setminus [0, \,+\infty)$.
We then have
\begin{gather}  \label{eqn:laplace-5-0}
\begin{split}
\Vert \varGamma_0(z) u \Vert_{H^{\mu, -s{\prime}}}
=
\Vert
\langle D \rangle^{\mu} \varGamma_0(z) u 
\Vert_{L^{2, -s^{\prime}}}                  
=
\Vert
 \varGamma_0(z) \langle D \rangle^{\mu}u 
\Vert_{L^{2, -s^{\prime}}}          \qquad   \\
\noalign{\vskip 4pt}
\le
\Vert
 \varGamma_0(z) 
\Vert_{B(0, s \, ; \,0, -s^{\prime})} 
\,
\Vert 
\langle D \rangle^{\mu}u
\Vert_{L^{2, s}}           \qquad\qquad \qquad \qquad           \\
\noalign{\vskip 4pt}   
=
\Vert
 \varGamma_0(z) 
\Vert_{B(0, s \, ; \,0, -s^{\prime})} 
\,
\Vert u \Vert_{H^{\mu, s}},  \qquad\qquad \qquad \qquad   \quad
\end{split}
\end{gather}
which implies (\ref{eqn:laplace-3}). 
In a similar fashion, 
we can prove (\ref{eqn:laplace-3+}).
\end{proof}

\begin{rem}  \label{rem:independence-1}
We should remark that ${H}^{\mu, \, s}({\mathbb R}^3)$ 
 in Lemma \ref{lem:lap-1} is
 a subset of ${L}^2({\mathbb R}^3)$ for $\mu \ge 0$, but
not necessarily for $\mu <0$. Thus, the domain
of $\widetilde\varGamma_0(z)$ 
 depends on $\mu$ and $s$. Nonetheless, we have
the unique representation of $\widetilde\varGamma_0(z)$ 
on  $S({\mathbb R}^3)$,
a dense subset of ${H}^{\mu, \, s}({\mathbb R}^3)${\rm :}
\begin{equation}   \label{eqn:kernelrep-1}
\widetilde\varGamma_0(z) u(x)
= 
\int_{{\mathbb R}^3} 
\frac{e^{i \sqrt{z} |x-y|} }{\, 4\pi |x-y| \,} u(y) \, dy,
\qquad
u \in S({\mathbb R}^3),
\end{equation}
for every $z \in \Pi_{(0, \, +\infty)}$,
where $\mbox{\rm Im}\sqrt{z} \ge 0$.
This representation, together with 
the fact that $S({\mathbb R}^3)$ is dense
in ${H}^{\mu, \, s}({\mathbb R}^3)$ for 
any pair of $\mu$ and $s$,
 ensures that  
the extension of 
$\widetilde\varGamma_0(z){\big|}_{S({\mathbb R}^3)}$ 
to
${H}^{\mu, \, s}({\mathbb R}^3)$ is independent of
 $\mu$ and $s$
in a certain sense.
However, we shall not discuss about the uniqueness
of the extension
any longer.
In the discussions below, we shall mostly
deal with the extension of 
$\widetilde\varGamma_0(z){\big|}_{S({\mathbb R}^3)}$ 
to
${H}^{-1, \, s}({\mathbb R}^3)$.
\end{rem}

\begin{lem} \label{lem:lap-1-1}
Let $s$, $s^{\prime}$ satisfy
the same assumption as in Lemma \ref{lem:hs},
 and let $\mu \in \mathbb R$.
 Then $\widetilde\varGamma_0(z)$ is a  
${B}(\mu -2 , s \, ; \,\mu, - s^{\prime})$-valued
continuous function
on $\Pi_{(0, \, +\infty)}$.
\end{lem}

\begin{proof} 
We first note that
\begin{equation}  \label{eqn:laplace-5}
\langle D \rangle^{2} \varGamma_0(z) u
= u + (z +1) \varGamma_0(z)  u
\end{equation}
for $u\in S({\mathbb R}^3)$ and 
$z\in \mathbb C \setminus [0, \,+\infty)$;
cf. Jensen and Kato \cite[Lemma 2.1]{JensenKato}.
(See (\ref{eqn:bracket-D}) for 
the definition of $\langle D \rangle$.) 
We then combine (\ref{eqn:laplace-5})
with Lemma \ref{lem:lap-1}, and obtain the conclusion
if we appeal to the fact that $S({\mathbb R}^3)$ is
dense in ${H}^{\mu -2, s}({\mathbb R}^3)$.
\end{proof}

\vspace{10pt}
 What we shall need in the rest of the paper
is a variant of 
Lemma \ref{lem:lap-1-1}, namely a version for
four-component
vector-valued functions, with $\mu = 1$ in the
form described in Proposition \ref{prop:lap-1-2} below. 
Thus
$\widetilde\varGamma_0(z)$ 
in Proposition \ref{prop:lap-1-2}  denotes
a copy of $\widetilde\varGamma_0(z)$; see (\ref{eqn:mapsto}).

\begin{prop} \label{prop:lap-1-2}
Let $s$, $s^{\prime}$ satisfy
the same assumption as in Lemma \ref{lem:hs}.
 Then $\widetilde\varGamma_0(z)$ is a  
${\mathcal B}(-1 , s \, ; \, 1, - s^{\prime})$-valued
continuous function
on $\Pi_{(0, \, +\infty)}$.
\end{prop}

\vspace{10pt}

\subsection{The resolvent of the free Dirac operator $H_0$}

In view of (\ref{eqn:reslvnDrc-0}),
it is convenient for us to introduce 
the following operator
valued-functions $\varOmega_0^{+}(z)$ 
defined on $\overline{{\mathbb C}}_{+}$
and $\varOmega_0^{-}(z)$
on $\overline{{\mathbb C}}_{-}$
as follows:
\begin{equation}       \label{eqn:laplace-6-0}
\varOmega_0^{\pm}(z) 
=\widetilde\varGamma_0(z^2), \quad z \in
 \overline{{\mathbb C}}_{\pm},
\end{equation}
in other words,
\begin{equation}   \label{eqn:laplace-6}
\varOmega_0^{\pm}(z) 
=
\begin{cases}  
\varGamma_0(z^2) 
& \mbox{if } \; z \in {\mathbb C}_{\pm}, 
\\ {}    & \\
\varGamma_0^{\pm}(\lambda^2) 
  & \mbox{if } \;  z= \lambda \ge 0,
\\ {}    & \\
\varGamma_0^{\mp}(\lambda^2) 
  & \mbox{if } \; z= \lambda \le 0.
\end{cases}
\end{equation}
It follows from Proposition \ref{prop:lap-1-2} 
 that $\varOmega_0^+(z)$ ( resp. $\varOmega_0^-(z)$) is
a $\mathcal{B}(-1 , s \, ; \,1, - s^{\prime})$-valued
continuous function on $\overline{\mathbb C}_+$
(resp. $\overline{\mathbb C}_-$). 
Also, it follows from (\ref{eqn:resolventLapc-3-1})
that
\begin{equation}   \label{eqn:addtional1}
\varOmega_0^+(0) = \varOmega_0^-(0)
=\widetilde\varGamma_0(0).
\end{equation}

In order to get expressions of the extended resolvents
of the free Dirac operator $H_0$ in terms of 
$\varOmega_0^{\pm}(z)$ introduced in (\ref{eqn:laplace-6}),
we shall exploit the decomposition
(\ref{eqn:reslvnDrc-0})
and a boundedness estimate of $H_0$ in some 
weighted Sobolev spaces which
is given as follow. 

\begin{lem}  \label{lem:dirac-bdd}
Let $\mu$ and $s^{\prime}$ be in $\mathbb R$. Then
\begin{equation*}   \label{eqn:dirac-bdd-2}
H_0 \in  \mathcal B (\mu, -s^{\prime} ; \mu-1, -s^{\prime}).
\end{equation*}
\end{lem}

\begin{proof}
To prove the lemma, it is sufficient to
show that
\begin{equation*} 
\langle x \rangle^{-s^{\prime}}
\langle D \rangle^{\mu -1}
D_j
\langle D \rangle^{-\mu}
\langle x \rangle^{s^{\prime}}
=   
\langle x \rangle^{-s^{\prime}}D_j
\langle D \rangle^{-1}
\langle x \rangle^{s^{\prime}}
\end{equation*}
($j=1$, $2$, $3$) is a bounded operator from 
$L^2({\mathbb R}^3)$ to 
$L^2({\mathbb R}^3)$.
This fact is a direct consequence of 
Umeda \cite[Lemma 2.1]{Umeda0}.
\end{proof}

\vspace{5pt}
\begin{lem}   \label{lem:lap-2}
Let $s$, $s^{\prime} > 1/2$, and $s + s^{\prime}>2$.
Then
$R_0(z) \in 
{\mathcal B}(-1, s \, ; \, 0, - s^{\prime})$ is
continuous in $z \in {\mathbb C}_{\pm}$.
Moreover, 
as ${\mathcal B}(-1, s \, ; \, 0, - s^{\prime})$-valued
functions,
 they
can possess continuous extensions 
$R_0^{\pm}(z)$ to
$\overline{\mathbb C}_{\pm}$ respectively,
and
\begin{equation} \label{eqn:dirac-1}
R_0^{\pm}(z) = (H_0 + z) \varOmega_0^{\pm}(z),  
 \quad z \in \overline{\mathbb C}_{\pm}. 
\end{equation}
\end{lem}

\begin{proof}
We shall give the proof only for 
$z \in \overline{\mathbb C}_{+}$.
The proof for $z \in \overline{\mathbb C}_{-}$ is
similar.

As was mentioned before Lemma \ref{lem:dirac-bdd}, 
$\varOmega_0^{+}(z)$
is
a $\mathcal{B}(-1 , s \, ; \,1, - s^{\prime})$-valued
continuous function on $\overline{\mathbb C}_+$.
Combining this fact with (\ref{eqn:reslvnDrc-0}),
 (\ref{eqn:reslvnDrc-1}), 
the definition (\ref{eqn:laplace-6-0}) 
(or (\ref{eqn:laplace-6}))
of $\varOmega_0^{+}(z)$, 
 Proposition \ref{prop:lap-1-2} and 
Lemma \ref{lem:dirac-bdd} with
$\mu=1$, 
we see that
$R_0(z)= (H_0 + z) \varOmega_0^{+}(z) \in 
{\mathcal B}(-1, s \, ; \, 0, - s^{\prime})$ for
any $z \in {\mathbb C}_{+}$. 
Now it is evident that 
the second assertion of the lemma follows  from
Proposition \ref{prop:lap-1-2} and 
Lemma \ref{lem:dirac-bdd} with
$\mu=1$.
\end{proof}

Combining (\ref{eqn:dirac-1}) with (\ref{eqn:addtional1}),
we obtain a corollary to Lemma \ref{lem:lap-2}.

\begin{cor}  \label{cor:origin0}
Under the same assumption and the same notation as in
Lemma \ref{lem:lap-2},
\begin{equation*}   \label{eqn:additional2}
R_0^+(0)=R_0^-(0) = H_0\widetilde\varGamma(0)
\quad
\mbox{in } \; 
{\mathcal B}(-1, s \, ; \, 0, - s^{\prime}).
\end{equation*}
\end{cor}

\vspace{4pt}

\begin{rem}
In \cite{IftMant}, Iftimovici and  M\u{a}ntoiu showed that
the limiting absorption principle
for the the free Dirac operator
$H_0= \alpha \cdot D + m \beta$, 
$m>0$, in ${\mathcal B}(0, 1 \, ; \, 0, - 1)$
holds  on the whole real line. With the result
exhibited in Lemma \ref{lem:lap-2},
together with the result in \cite{IftMant},
the limiting absorption principle
for the the free Dirac operator
$H_0= \alpha \cdot D + m \beta$
has been established 
for all $m \ge 0$.
\end{rem}

\begin{lem} \label{lem:ker3-2}
For $f \in { \mathcal S }$ and $z \in {\mathbb C}_{\pm}$
\begin{eqnarray} \label{eqn:3-2}
{}&R_0(z)f (x) \qquad \qquad \qquad \qquad
  \qquad \qquad \qquad \qquad \qquad \qquad \\
{}&= \displaystyle{
\int_{{\mathbb R}^3}  \Big(
i\, \frac{\alpha \cdot (x-y)}{|x-y|^2} \pm 
z \, \frac{\alpha \cdot (x-y)}{|x-y|}  + zI_4 \Big)
\frac{e^{\pm i z |x-y|}}{4\pi |x-y|}
} f(y) \, dy.   \nonumber
\end{eqnarray}
\end{lem}
\begin{proof}
We first recall (\ref{eqn:kernelrep-1}),
which we can write as
\begin{equation}   \label{eqn:kernelrep-1-1}
\varGamma_0(z^2) f(x)
= 
\int_{{\mathbb R}^3} 
\frac{e^{ \pm i z |y|} }{\, 4\pi |y| \,} f(x-y)  \, dy,
\qquad
 f \in {\mathcal S}, \;\; z \in {\mathbb C}_{\pm}.
\end{equation}
We then combine (\ref{eqn:laplace-6-0}) and 
(\ref{eqn:dirac-1}), and make differentiation
under the integral sign in (\ref{eqn:kernelrep-1-1}),
which gives
\begin{gather}   \label{eqn:3-2-0}
\begin{split}
{}&R_0(z)f (x) = 
\int_{{\mathbb R}^3}  
\frac{e^{\pm i z |y|}}{4\pi |y|}
             (\alpha \cdot D_x +z I_4)      f(x-y) \, dy ,
\quad z  \in {\mathbb
C}_{\pm}.  
\end{split} 
\end{gather}
Noting the fact that 
\begin{equation*}    \label{eqn:differentiation}
D_x f(x-y)= - D_y\big( f(x-y) \big),
\end{equation*}
and making integration by parts on the right hand
of (\ref{eqn:3-2-0}) 
implies that
\begin{gather}   \label{eqn:3-2-1}
\begin{split}
R_0(z)f (x)  &= 
\int_{{\mathbb R}^3}  \Big(
i\, \frac{\alpha \cdot y}{|y|^2} \pm 
z \, \frac{\alpha \cdot y}{|y|}  + zI_4 \Big)
\frac{e^{\pm i z |y|}}{4\pi |y|}
                   f(x-y) \, dy      \\
\noalign{\vskip 4pt}
{}& \qquad\quad (z  \in {\mathbb C}_{\pm}).  
\end{split} 
\end{gather}
A change of variables in (\ref{eqn:3-2-1})
yields (\ref{eqn:3-2}).
(See also Thaller \cite[p.39]{Thaller}.)
\end{proof}

\begin{prop} \label{prop:riesz3-3}
For $f \in { \mathcal S }$ 
\begin{equation*} \label{eqn:3-3}
R_0^{+}(0) f = R_0^{-}(0) f = A f,
\end{equation*}
where  $A$ is the singular integral operator 
defined by {\rm(\ref{eqn:sio-1})}.
\end{prop}
\begin{proof}
In view of Corollary \ref{cor:origin0},
we only need to give the proof for $R_0^+(0)$. 

Let $f \in { \mathcal S }$, and 
let $\{ z_n \} \subset {\mathbb C}_{+}$
be a sequence such that $z_n \to 0$ as $n \to \infty$.
It follows from Lemma \ref{lem:lap-2} that
$R_0(z_n)f \to R^+_0(0)f$ in 
${\mathcal L}^{2, -s^{\prime}}$ as $n \to
\infty$.
This fact implies that there exists a subsequence
$\{ z_{n^{\prime}} \} \subset\{ z_n \}$ 
such that
\begin{equation}   \label{eqn:convergence-1}
R_0(z_{n^{\prime}})f(x) \to R^+_0(0)f(x)
\;\; \mbox{ a.e. } x \in {\mathbb R}^3 
\;\;\;  
\mbox{ as } n^{\prime} \to \infty.
\end{equation}

On the other hand, Lemma \ref{lem:ker3-2},
together with  Lebesgue's
 convergence theorem, implies
that 
\begin{equation}     \label{eqn:convergence-2}
R_0(z_{n})f(x) \to 
\int_{{\mathbb R}^3}  
i\, \frac{\alpha \cdot (x-y)}{4\pi |x-y|^3}
 f(y) \, dy = Af(x)
\;\;\;  
\mbox{ as } n  \to \infty
\end{equation}
for each $x \in {\mathbb R}^3$.
The conclusion of the proposition now 
follows from (\ref{eqn:convergence-1}) and 
(\ref{eqn:convergence-2}).
\end{proof}

\begin{lem} \label{lem:riesz3-4}
Let $s$, $s^{\prime} > 1/2$, and $s + s^{\prime}>2$.
 Then
$A$ can be continuously extended to
an operator in 
${\mathcal B}(-1, s \, ; \, 0, - s^{\prime})$. 
\end{lem}
\begin{proof}
Since $\mathcal S$ is dense in
${\mathcal H}^{-1, s}$, Lemma \ref{lem:lap-2}
and Proposition \ref{prop:riesz3-3} directly
imply the lemma.
\end{proof}

In the rest of the paper, we shall denote the extension
in Lemma \ref{lem:riesz3-4} 
by $A$ again. Thus we have
\begin{equation*}
R_0^+(0) =R_0^-(0) = A \; 
\mbox{ in } \; 
{\mathcal B}(-1, s \, ; \, 0, - s^{\prime}).
\end{equation*}

\begin{prop} \label{prop:riesz3-5-0}
Let $s > 1/2$. Then 
\begin{equation}   \label{eqn:adjoint-1}
H_0 A g = g
\end{equation}
for all $g \in {\mathcal L}^{2, s}$.
\end{prop}

\begin{proof}
Let $g \in {\mathcal L}^{2, s}$ be given.
We  then start with the fact that
\begin{equation}   \label{eqn:adjoint-1+}
(H_0 - z) R_0(z)g =g  \quad ( \forall z \in {\mathbb C}_+ ).
\end{equation}
Choose $s^{\prime} > 1/2$ so that
$s + s^{\prime} >2$. 
We see from Lemmas \ref{lem:lap-2}, \ref{lem:riesz3-4}
and Proposition \ref{prop:riesz3-3} 
that
\begin{equation}    \label{eqn:pair-4}
R_0(i/n) g \to R_0^+(0)g = Ag  \;\; \mbox{ \rm in }  \;\; 
{\mathcal L}^{2, -s^{\prime}}
\quad
{\mbox \rm as } \;\; n   \to \infty .
\end{equation}
Lemma \ref{lem:dirac-bdd}, with $\mu=0$,
and (\ref{eqn:pair-4})
imply that
\begin{equation}    \label{eqn:pair-5}
\Big( H_0 -\frac{i}{\, n \,} \Big) 
 R_0\Big( \frac{i}{\, n \,} \Big) g 
 \to  H_0 Ag  \;\; \mbox{ \rm
in }  \;\;  {\mathcal H}^{-1, -s^{\prime}}
\quad
{\mbox \rm as } \;\; n \to \infty.
\end{equation}
Since, by (\ref{eqn:adjoint-1+}),
\begin{equation*}    \label{eqn:pair-5+}
\Big( H_0 - \frac{i}{\, n \,} \Big) 
 R_0 \Big( \frac{i}{\, n \,}\Big) g 
=g  \quad  \mbox{ for }\forall n,
\end{equation*}
we find that 
(\ref{eqn:pair-5}) yields (\ref{eqn:adjoint-1}).
\end{proof}

\vspace{5pt}

We shall need  
Lemma 2.4 of Jensen and Kato \cite{JensenKato},
which we shall rewrite in a suitable form
to our setting (cf. Lemma \ref{lem:riez3-5-7} below),
where the operators
$-\Delta$
and 
$\widetilde\varGamma_0(0)$
act on four-component vector functions. 
The reader should note
that $\widetilde\varGamma_0(0)$
is the same as $G_0$ in 
Jensen-Kato's paper. 
See (\ref{eqn:resolventLapc-1}) and
(\ref{eqn:resolventLapc-3-1}).

\begin{lem}[\textbf{Jensen-Kato}]   \label{lem:riez3-5-7}
Let $s >1/2$. Then

\vspace{3pt}

{\rm(i)} $(-\Delta) \widetilde\varGamma_0(0) g = g$ \ 
 for all \  $g \in {\mathcal H}^{-1, s}.$

\vspace{2pt}

{\rm(ii)} $ \widetilde\varGamma_0(0) (-\Delta)f = f \;$ 
if $\; f \in  {\mathcal L}^{2, -3/2} \ and \ 
(-\Delta)f \in {\mathcal H}^{-1, s}$.
\end{lem}

\vspace{4pt}

\begin{prop} \label{prop:riez3-5-8}
Let $s >1/2$. Then \ 
 $\widetilde\varGamma_0(0) H_0 \, g = Ag$  \ 
 for all $g \in {\mathcal L}^{2, s}.$
\end{prop}
\begin{proof}
Let $g \in {\mathcal L}^{2, s}$ be given.
Noting that  $H_0^2 = -\Delta$ (cf. (\ref{eqn:diracsquare})),
 we have
\begin{equation}   \label{eqn:identity-10}
(-\Delta) A g 
=H_0 (H_0 Ag)            
=H_0 \, g  
\end{equation}
where we have used Proposition \ref{prop:riesz3-5-0}
in the second equality.
Since $H_0 \, g \in {\mathcal H}^{-1, s}$ by Lemma \ref{lem:dirac-bdd},
it follows from (\ref{eqn:identity-10}) that
$(-\Delta) A g  \in {\mathcal H}^{-1, s}$.

On the other hand, we find, by Lemma \ref{lem:riesz3-4},
that $Ag \in {\mathcal L}^{2, -3/2}$,
because we can choose $s^{\prime}$ so 
that $1/2 < s^{\prime} \le 3/2$ \ and
$s + s^{\prime} >2$. (Choose 
$s^{\prime}$ so 
that $\max (s, \, 2-s) < s^{\prime} \le 3/2$.)

Now we can apply Lemma \ref{lem:riez3-5-7}(ii) with $f$ 
replaced by $Ag$, and obtain
\begin{equation}   \label{eqn:identity-11}
\widetilde\varGamma_0(0)(-\Delta) A g 
=Ag. 
\end{equation}
It follows from (\ref{eqn:identity-10}) that
the left hand side of (\ref{eqn:identity-11})
 equals $\widetilde\varGamma_0(0)
H_0g$. This proves the conclusion of the proposition.
\end{proof}
\vspace{2pt}
\begin{proof}[{\bf Proof of Theorem \ref{thm:riesz3-5}}]
Put
\begin{equation*}  \label{eqn:identity-12}
g= H_0 f.
\end{equation*}
By assumption of the theorem, we see that $g \in {\mathcal L}^{2, s}$ 
for some $\, s > 1/2$.
It follows from Proposition \ref{prop:riez3-5-8}
that $Ag = \widetilde\varGamma_0(0) H_0 \, g$, {\it i.e.},
\begin{equation*}    \label{eqn:identity-13}
AH_0f = \widetilde\varGamma_0(0) H_0 H_0f
  = \widetilde\varGamma_0(0) (-\Delta)f.
\end{equation*}
Since $(-\Delta)f=H_0g \in {\mathcal H}^{-1, s}$  by
Lemma \ref{lem:dirac-bdd}, 
it follows from assertion (ii) of 
Lemma \ref{lem:riez3-5-7} that
$\widetilde\varGamma_0(0) (-\Delta)f=f$. Thus $AH_0f =f$.
\end{proof}

\vspace{15pt}

\section{Proof of the main theorems}  \label{sec:pfthm}

\begin{proof}[{\bf Proof of Theorem \ref{thm:th-2-2}}]
We first prove assertion (i).
Let $f$ be a zero mode of the operator $H$ in
(\ref{eqn:1-1}). Then we have
\begin{equation}   \label{eqn:proof-1}
Hf = \big( \alpha \cdot D + Q(x) \big) f =0,
\quad
f \in \mbox{Dom}(H) = {\mathcal H}^1 .
\end{equation}
It follows from (\ref{eqn:proof-1}) and Assumption (A)
that
\begin{equation}   \label{eqn:proof-2}
H_0f =(\alpha \cdot D) f = - Q(x) f
 \in 
{\mathcal L}^{2, \rho}.
\end{equation}
(Recall  (\ref{eqn:reslvnDrc-00}) for the definition of $H_0$.)
Since $\rho >1 >1/2$ by assumption of the
theorem, we can apply Theorem \ref{thm:riesz3-5} 
to (\ref{eqn:proof-2})
and get
\begin{equation}   \label{eqn:proof-3}
f = AH_0 f = - A  Q f.
\end{equation}
It follows from (\ref{eqn:proof-3}) and
 Lemma \ref{lem:siolem-1} 
that
$f \in {\mathcal L^2} \cap {\mathcal L}^6$.
It follows from (\ref{eqn:proof-3}) 
again and Lemma \ref{lem:siolem-3} 
that $f \in {\mathcal L}^{\infty}$.
This fact, together with  (\ref{eqn:proof-3})
and Assumption (A),
implies that
\begin{align}   
| f(x) |
\le
\frac{3}{4\pi}
 \int_{{\mathbb R}^3}
 \frac{1}{\, |x-y|^2 \,} \,
        | Q(y) f(y) | \, dy  \label{eqn:proof-4-1}  \\
 \le 
C \Vert f \Vert_{{\mathcal L}^{\infty}} 
\int_{{\mathbb R}^3}
 \frac{1}{\, |x-y|^2 \, \langle y \rangle^{\rho} \,} 
\, dy.                         \label{eqn:proof-4}
\end{align}
Noting that $\rho >1$ by assumption, 
and  applying  the inequalities in (\ref{eqn:ekku})
 to the integral 
in (\ref{eqn:proof-4}),
we get
\begin{equation}   \label{eqn:pe-1}
|f(x)|
\le
C \Vert f \Vert_{{\mathcal L}^{\infty}} 
\begin{cases}  
\langle x \rangle^{-\rho + 1} 
& \mbox{if } \; 1  < \rho < 3, 
\\ {}    & \\
\langle x \rangle^{-2} \log ( 1 + \langle x \rangle )
  & \mbox{if } \; \rho = 3,
\\ {}    & \\
\langle x \rangle^{-2}
  & \mbox{if } \;  \rho > 3.
\end{cases}
\end{equation}
If $\rho >3$, we have already obtained the desired estimate.
If $1 < \rho \le 3$, we plug
 the inequalities in (\ref{eqn:pe-1})
into (\ref{eqn:proof-4-1}).
We thus get
\begin{equation}   \label{eqn:pe-2}
|f(x)|
\le
C \Vert f \Vert_{{\mathcal L}^{\infty}} 
\int_{{\mathbb R}^3}
 \frac{1}{\, |x-y|^2 \, \langle y \rangle^{2\rho -1} \,} 
\, dy  
\end{equation}
if $1  < \rho < 3 $,
and
\begin{equation}    \label{eqn:pe-3} 
|f(x)|
\le
C \Vert f \Vert_{{\mathcal L}^{\infty}} 
\int_{{\mathbb R}^3}
 \frac{\log(1+ \langle y \rangle )}
{\, |x-y|^2 \, \langle y \rangle^{\rho+2} \,} 
\, dy  
\end{equation}
if $\rho = 3$.
We find that the inequalities in (\ref{eqn:ekku})
 applied to the
integrals  in (\ref{eqn:pe-2}) and (\ref{eqn:pe-3}) 
yields 
\begin{equation}   \label{eqn:pe-4}
|f(x)|
\le
C \Vert f \Vert_{{\mathcal L}^{\infty}} 
\begin{cases}  
\langle x \rangle^{-2(\rho - 1)} 
& \mbox{if } \; 1  < \rho < 2, 
\\ {}    & \\
\langle x \rangle^{-2} \log ( 1 + \langle x \rangle )
  & \mbox{if } \; \rho = 2,
\\ {}    & \\
\langle x \rangle^{-2}
  & \mbox{if } \;   2 < \rho \le 3.
\end{cases}
\end{equation}
Hence, if  $ 2 < \rho \le 3$, we have shown the 
desired estimate.  If $ 1 < \rho \le 2$, we repeat
the same argument again, actually as many times as 
we wish.
Summing up, we  can obtain the estimate
\begin{equation}   \label{eqn:pe-5}
|f(x)|
\le
C_N \Vert f \Vert_{{\mathcal L}^{\infty}} 
\begin{cases}  
\langle x \rangle^{-N(\rho - 1)} 
& \mbox{if } \; 1  < \rho <  1 + \displaystyle{\frac{2}{N}}, 
\\ {}    & \\
\langle x \rangle^{-2} \log ( 1 + \langle x \rangle )
  & \mbox{if } \; \rho =  1 + \displaystyle{\frac{2}{N}},
\\ {}    & \\
\langle x \rangle^{-2}
  & \mbox{if } \;   1 + \displaystyle{\frac{2}{N}} < \rho .
\end{cases}
\end{equation}
for any positive integer $N$, where $C_N$ is a constant
depending on $N$. It is straightforward that 
for a given $\rho >1$ in Assumption
(A), we can choose $N$ so that 
$1 + (2/N) < \rho$.
This fact, together with (\ref{eqn:pe-5}),
 implies assertion (i).

We next prove assertion (ii) by
utilizing  (\ref{eqn:proof-3}):
\begin{equation*}   \label{eqn:proof-4+1}
f(x) = 
-\displaystyle{
\int_{{\mathbb R}^3}
  i \, \frac{\alpha \cdot (x - y)}{4\pi |x-y|^3} Q(y) f(y) \, dy.
}
\end{equation*}
Let $x_0$ be any point in ${\mathbb R}^3$, and
let
$\varepsilon >0$ be given.
We choose  $r>0$ so that
\begin{equation}   \label{eqn:proof-4+2}
\displaystyle{
\int_{|y| \le 2r}
  \frac{1}{\, |y|^2 \,} \, dy < \varepsilon}.
\end{equation}
We then decompose $f(x)$ into two
parts:
\begin{gather}   \label{eqn:proof-4+3}
\begin{split}
f(x)
&= - \Big( 
\int_{B(x, \,2r)} 
+ \int_{E(x,\,2r)}  
\Big)
 \, i \, 
   \frac{\alpha \cdot (x - y)}{4\pi |x-y|^3} Q(y) f(y) \, dy   \\
\noalign{\vskip 4pt}
&=: f_b(x) + f_e(x),     
\end{split}
\end{gather}
where
\begin{equation*}   \label{eqn:proof-4+4}
B(x, \, 2r) = \{ \, y \; \big| 
\,|x - y| \le 2r \, \}, 
\quad 
E(x, \,2r) =  \{ \, y \; \big| 
\,|x - y| > 2r \, \}.
\end{equation*}
Since each $\alpha_j$ is a unitary matrix,
it follows from (\ref{eqn:proof-4+2})
and (\ref{eqn:proof-4+3}) that
\begin{equation}    \label{eqn:proof-4+5}
|f_b(x)| < \frac{3}{\, 4\pi \,} C_q C_f \,\varepsilon  
\quad \mbox{ for }\; \forall x
\in {\mathbb R}^3 ,
\end{equation}
where $C_q$ is a constant determined 
by  (\ref{eqn:2-1}) in Assumption (A)
and $C_f$ is a constant described in 
the inequality (\ref{eqn:2-2}),
which we have just proved in
the first half of the proof.
It follows from the 
definition of $f_e(x)$ that
\begin{gather}   \label{eqn:proof-4+6}
\begin{split}
{}&f_e(x) - f_e(x_0)    \\
\noalign{\vskip 4pt}
&= 
\int_{{\mathbb R}^3} \!
\Big\{
1_{E(x_0, \, 2r)}(y) 
\frac{\alpha \cdot (x_0 - y)}{4\pi |x_0-y|^3}  
-
1_{E(x, \, 2r)}(y) 
\frac{\alpha \cdot (x - y)}{4\pi |x-y|^3} 
\Big\}   \\
\noalign{\vskip 4pt}
&\qquad\qquad \qquad\qquad \qquad\qquad \qquad\qquad \qquad\qquad 
   \times Q(y) f(y) \, dy .
\end{split}
\end{gather}
To apply Lebesgue's  convergence theorem
to the integral in (\ref{eqn:proof-4+6}),
we need the following two facts that
\begin{equation}  \label{eqn:proof-4+7}
|x_0 - y | \ge  r
\quad \mbox{ if } \; 
 |x_0 - x| <r, \; |x-y| > 2r.
\end{equation}
and that
\begin{equation}  \label{eqn:proof-4+8}
|x - y | \ge  \frac{2}{\, 3 \,} |x_0 -y|
\quad \mbox{ if } \; 
 |x_0 - x| <r, \; |x-y| > 2r
\end{equation}
(use the inequality $|x-y| \ge |x_0 -y| - |x_0 -x|$).
We can deduce from
(\ref{eqn:proof-4+7}) and (\ref{eqn:proof-4+8}) that
\begin{gather}    \label{eqn:proof-4+9}
\begin{split}
{}&\Big| 
1_{E(x, \, 2r)}(y) 
\frac{\alpha \cdot (x - y)}{4\pi |x-y|^3} Q(y) f(y)  \Big| \\  
&\quad \le
1_{E(x_0, \, r)}(y) 
\frac{3}{\, 4\pi \,} 
\Big(
\frac{2}{\, 3 \,} |x_0 -y|
\Big)^{-2}   \,
 \big|  Q(y) f(y)  \big|
\end{split}
\end{gather}
whenever $|x_0 - x| < r$.
It is
straightforward that
the estimate (\ref{eqn:proof-4+9}) 
implies 
\begin{gather}   \label{eqn:proof-4+10}
\begin{split}
{}&\big|
\mbox{the integrand in (\ref{eqn:proof-4+6})}
\big|  \\
\noalign{\vskip 4pt}
&\le 
1_{E(x_0, \, r)}(y) 
\frac{3}{\, 4\pi \,} 
\Big(
1+ (\frac{3}{\, 2 \,})^2 
\Big)
|x_0 -y|^{-2}   \,
 \big|  Q(y) f(y)  \big| 
\end{split}
\end{gather}
whenever $|x_0 - x| < r$.
In view of
(\ref{eqn:2-1}) in Assumption (A)
and the inequality (\ref{eqn:2-2}),
the function on the right hand side
of (\ref{eqn:proof-4+10})
is integrable on ${\mathbb R}^3$.
Thus, we can apply  
 Lebesgue's  convergence theorem
to the integral in (\ref{eqn:proof-4+6}),
and conclude that
\begin{equation}    \label{eqn:proof-4+11}
\lim_{x \to x_0} \big( f_e(x) - f_e(x_0) \big)=0.
\end{equation}
Combining (\ref{eqn:proof-4+11})
with 
both (\ref{eqn:proof-4+3}) and (\ref{eqn:proof-4+5}) yields
\begin{equation*}   \label{eqn:proof-4+12}
 \limsup_{x \to x_0}
\big| f(x) - f(x_0) \big| 
\le 
2\times
\frac{3}{\, 4\pi \,} C_q C_f \,\varepsilon .     
\end{equation*}
Since $\varepsilon$ was arbitrary,
this completes the proof of assertion (ii).
\end{proof}

\vspace{15pt}

\begin{proof}[{\bf Proof of Theorem \ref{thm:th-2-1}}]
Let $f$ satisfy the assumption of the theorem:
 $f\in {\mathcal L}^{2,-s}$
for  some $s$ with $0 < s  \le \min\{3/2, \, \rho -1\}$. 
In the same manner as in (\ref{eqn:proof-2}) and 
(\ref{eqn:proof-3}), we can show that
\begin{equation}  \label{eqn:proof-7}
H_0f=(\alpha \cdot D) f =- Qf \in  
{\mathcal L}^{2,\rho-s},
\end{equation}
and that
\begin{equation}  \label{eqn:proof-8}
f= -A Qf.
\end{equation}
Note that $s \le 3/2$ and $\rho -s \ge 1 >1/2$, which we have used 
to apply Theorem \ref{thm:riesz3-5}
 in showing 
 (\ref{eqn:proof-8}).
Since $Qf \in  
{\mathcal L}^{2,\rho-s}$, $\rho-s \ge 1$,
we see from  (\ref{eqn:proof-8}) and 
Lemma \ref{lem:siolem-2}
that $f \in {\mathcal L}^2$.
This fact, together with (\ref{eqn:proof-7})
and Lemma \ref{lem:sobolev-1}, gives 
the conclusion of the theorem.
\end{proof}

\vspace{20pt}

{\bf Acknowledgment:}  T.U. would like to express his gratitude to 
Michael Loss
for the hospitality during his visit to 
Georgia Institute of Technology, USA, in
April, 2002. Discussions with Michael were one of the 
motivations of the present
paper.
Also, he would like to express his thanks to the Department
of Mathematics, the University of 
Alabama at Birmingham, USA, 
for their hospitality.
Part of the present paper was done during
his stay there in March and September, 2006.
Finally the authors
 appreciate  invaluable comments by Michael Loss, Kenji
 Yajima and the referee. Kenji's comments helped us  improve
the main theorems of the previous version of 
the present paper.


\end{document}